\documentclass[10pt]{amsart}
\usepackage{graphicx}
\usepackage{color}\usepackage{amsmath,amsthm,amssymb,amsfonts}
\usepackage[english]{babel}
\usepackage{amsmath}
\usepackage{amssymb}
\newtheorem{teo}{Theorem}[section]

\newtheorem{cor}[teo]{Corollary}

\newtheorem{prop}[teo]{Proposition}
\theoremstyle{definition}

\theoremstyle{remark}

\numberwithin{equation}{section}

\newcommand{\N}{\mathbb{N}}

\newcommand{\C}{\mathbb{C}}
\newcommand{\D}{\mathbb{D}}
\newcommand{\B}{\mathcal{B}}
\newcommand{\dem}{\noindent{\textsf{Proof.} }}
\newcommand{\hinf}{H^{\infty}}
\newcommand{\linf}{\ell_\infty}
\newcommand{\Bpre}{^*\B_v^\infty}
\newcommand{\Bpreo}{^*\widetilde{\B}_v^\infty}
\newcommand{\Binf}{\B_{v}^\infty}
\newcommand{\Bo}{\B_v^0}
\newcommand{\Binfo}{\widetilde{\B}_{v}^\infty}
\newcommand{\Boo}{\widetilde{\B}_v^0}
\renewcommand{\phi}{\varphi}
\newcommand{\eps}{\varepsilon}

\title{On interpolating sequences for Bloch type spaces}
\author{Alejandro Miralles}
\author{Mario P. Maletzki}

\address{Departament de Matemàtiques and IMAC, Universitat Jaume I, Castelló (Spain). \emph{e}.mail: mirallea@uji.es and al411086@uji.es}

\thanks{The first author was supported by Project PGC2018-094431-B-I00 (MICINN. Spain)}


\begin{document}

\maketitle	

\begin{abstract}
When we deal with $H^{\infty}$, it is known that $c_0-$interpolating sequences are interpolating and it is sufficient to interpolate idempotents of $\linf$ in order to interpolate the whole $\linf$. We will extend these results to the frame of interpolating sequences for Bloch type spaces $\Binf$ and study the connection between the interpolating operators on $\Binf$ and $\Bo$. Furthermore, for some particular weights $v$, we will provide examples of interpolating sequences for $\Binf$ whose constant of separation is as close to 0 as desired. 
%
	\end{abstract}

\section{Introduction and Background}

Let $\D$ be the open unit disk of the complex plane $\C$ and $H(\D)$ the space of complex analytic functions on $\D$. In this paper we will investigate sequences $(z_n) \subset \D$ which are interpolating for the derivative of functions in Bloch type spaces (see \cite{At}, \cite{BGLM}). It is also possible to study these sequences for Bloch type spaces that do not take into account the derivative of the function. For classical Bloch spaces, this has been done in \cite{BN}. \medskip

Let $v$ be a weight, that is, a strictly positive continuous function on $\D$ and suppose that $v$ is typical: $v$ is radial ($v(z)=v(|z|)$ for any $z \in \D$), non-increasing and $\lim_{|z| \rightarrow 1} v(z)=0$. \medskip

The Bloch type spaces $\B_v^{\infty}$ and $\B_v^0$ are defined by:
$$\B_v^{\infty}=\{ f \in H(\D) : \|f\|_{\B_v^{\infty}}:=|f(0)|+\sup_{z \in \D} v(z) |f'(z)| < +\infty \}$$
$$\B_v^0=\{f \in \B_v^{\infty} : \lim_{|z| \rightarrow 1^{-}} v(z)|f'(z)|=0 \}.$$

\noindent It is clear that $\Bo$ is a closed subspace of $\Binf$. We will also consider $\Binfo$ and $\Boo$, the closed subspaces of $\Binf$ and $\Bo$ respectively consisting in functions $f$ satisfying $f(0)=0$. It is also clear that $\Boo$ is a closed subspace of $\Binfo$. For typical weights $v$, it is well-known that the closed unit ball of $\Boo$ is dense with respect to the compact-open topology ($co-$topology) in the closed unit ball of $\Binfo$. Indeed, for $f$ in the closed unit ball of $\Binfo$, the functions $f_r(z)=f(rz)$ belong to the closed unit ball of $\Boo$ for $0 < r < 1$ and $f_r \stackrel{co}{\longrightarrow} f$. \medskip

Recall that $\hinf$ is the classical space of bounded analytic functions $f: \D \to \C$ endowed with the supremum norm $\| \cdot \|_{\infty}$. If $v(z)=1-|z|^2$, then $\B_{\infty}^v$ and $\B_{0}^v$ become the classical Bloch space $\B$ and the little Bloch space $\B_0$ respectively. It is well-known that $\hinf$ is properly contained in $\B$. This will remain true if we deal with weights $\mathcal{O}(1-|z|^2)$, but it is not true in general. Take $v_{\alpha}(z)=(1-|z|)^{\alpha}$ for $0 < \alpha < 1$ and $f(z)=(1-z)^\beta$ for $0 < \beta < 1-\alpha$. Then $f$ belongs to $\hinf$ but $|f'(z)|v_{\alpha}(z) \rightarrow \infty$ when $z \rightarrow 1$, so $f \notin \Binf$. \medskip

Let us first adapt some results due to Bierstedt and Summers done for the weighted Banach spaces of analytic functions $H_v^\infty$ and $H_v^0$ to the frame of Bloch type spaces \cite{BS}. Denote by $i: \Binfo \to (\Binfo)^{**}$ the natural inclusion of $\Binfo$ into its bidual $(\Binfo)^{**}$ given by $i(f)(u)=u(f)$ for $f \in \Binfo$ and $u \in (\B_v^\infty)^{*}$. In \cite{EGLN} it has been pointed out that the closed unit ball of $\Binfo$ is $co-$compact, so using the Dixmier-Ng Theorem \cite{N} we obtain that the space:
$$\Bpreo=\{ \ell \in (\Binfo)^{*} : \ell|_{B_{\Binfo}} \mbox{ is $co-$continuous} \}$$
\noindent endowed with the norm induced by the dual space $(\Binfo)^*$ is a Banach space and the map $f \in \Binfo \mapsto i(f)|_{\Bpreo} \in (^*\Binfo)^*$ is an onto isometric isomorphism. In particular, $\Bpreo$ is a predual of $\Binfo$. \medskip

Consider the evaluation functionals $\delta'_z$ given by $\delta'_z(f)=f'(z)$ for any $z \in \D$, which clearly belong to $\Bpreo$. By the Hahn-Banach theorem it follows that the linear span of $\{ \delta'_z : z \in \D \}$ is norm dense in $^*\B_v^\infty$. Therefore, $^*\B_v^\infty$ is separable since it is sufficient to consider evaluations $\delta'_z$ for $z=p+iq$, where $p,q \in \mathbb{Q}$. \medskip

Following the argument given by Bierstedt and Summers \cite{BS} we show:

\begin{prop}
The space $\Bpreo$ is isometrically isomoprhic to $(\Boo)^*$ and the restriction map $\widetilde{R}:$ $\Bpreo  \to (\Boo)^* $ given by $\widetilde{R}(\ell)=\ell|_{\Boo}$ is an onto isometric isomorphism. In particular, $(\Boo)^{**}$ is isometrically isomorphic to $\Binfo$.
\end{prop} 

\dem The map $\widetilde{R}$ is well-defined  since for any $\ell \in$ $\Bpreo \subset (\Binfo)^*$ it follows that $\|\widetilde{R}(\ell)\| \leq \| \ell\|$. First we prove that $\widetilde{R}$ is surjective. Consider $\ell \in (\Boo)^*$. Notice that $\Boo$ is isometrically isomorphic to a closed subspace of $C_0(\D)$, the space of continuous functions on the closed unit disk which vanish in the boundary, via $f \mapsto v f'$. By the Hahn-Banach theorem and the Riesz representation theorem, there is a bounded Radon measure $\mu$ on $\D$ such that:
$$\ell(f)=\int_\D v f' d \mu \ \mbox{ for } f \in \Boo.$$

\noindent Define $\tilde{\ell}(f)=\int_\D v f' d\mu$ for all $f \in \Binfo$ which is clearly well-defined and satisfies $\tilde{\ell}|_{\Boo}=\ell$. It follows from the Lebesgue bounded convergence theorem that $\tilde{\ell}|_{\Binfo}$ is $co-$continuous, so $\widetilde{R}$ is surjective. Since the closed unit ball of $\Boo$ is $co-$dense in the closed unit ball of $\Binfo$, we conclude that $\widetilde{R}$ is an isometry. \qed \bigskip

\begin{cor} \label{rem12}
The space $(\Bo)^*$ is isometrically isomorphic to $\Bpre$ and $(\Bpre)^*$ is isometrically isomorphic to $\Binf$. In particular, the space $(\Bo)^{**}$ is isometrically isomorphic to $\Binf$.
\end{cor}
  
  \dem Notice that $\Bo$ is isometrically isomorphic to $(\Boo \times \C,\| \cdot\|_1)$, so $(\Bo)^*$ is isometrically isomorphic to $\Bpre:=(\Bpreo \times \C,\| \cdot\|_\infty)$. The dual of this space is isometrically isomorphic to $(\Binfo \times \C,\| \cdot\|_1)$ which in turn is isometrically isomorphic to $\Binf$ and we conclude that $(\Bo)^{**}$ is isometrically isomorphic to $\Binf$. \qed \bigskip
  
\section{Interpolating sequences for Bloch type spaces}
Recall that the pseudohyperbolic distance for $z,w \in \D$ is given by:
$$\rho(z,w)=\left| \frac{z-w}{1-\bar{z}w}\right|.$$
A sequence $(z_n) \subset \D$ is said to be separated if there exists $\delta >0$ such that:
\begin{eqnarray} \label{sep}
\rho(z_n,z_k) \geq\delta \ \mbox{ for any } n \neq k,\end{eqnarray}

\noindent and we define its constant of separation as $r:=\inf_{n\neq k}\rho(z_n,z_k)$. \medskip

\noindent A sequence $(z_n) \subset \D$ is said to be interpolating for $\hinf$ if for any $(a_n) \in \ell_{\infty}$ there exists $f \in \hinf$ such that $f(z_n)=a_n$ for any $n \in \N$. The most important result on interpolating sequences for $\hinf$ is the classical Carleson's Theorem \cite{Carleson}, which states that $(z_n) \subset \D$ is interpolating for $\hinf$ if and only if $(z_n)$ is uniformly separated, that is, if there exists $\delta > 0$ such that $\inf_{k\in\mathbb{N}} \prod_{n \neq k} \rho(z_n,z_k) \geq\delta$. \medskip

A sequence $(z_n)\subset\mathbf{D}$ is said to be interpolating for the Bloch type space $\Binf$ if for any $(a_n) \in \linf$ there exists $f\in \Binf$ such that $v(z_n)f'(z_n)=a_n \mbox{ for any } n \in \N.$ 
We define the interpolating operator $T:\Binf \to \ell^\infty$ by $T(f)=(v(z_n)f'(z_n))$, which is clearly well-defined and linear. Notice that $(z_n)$ is interpolating for $\Binf$ if and only $T$ is surjective. If $(z_n)\subset\mathbf{D}$ satisfies $|z_n| \rightarrow 1$, then the interpolating operator $T|_{\Bo}$ maps $\Bo$ into $c_0$ since $f'(z_n)v(z_n) \rightarrow 0$ when $n \rightarrow \infty$. \medskip

For $\hinf$ and Bloch type spaces, we can also consider $c_0-$interpolating sequences just by considering sequences $(a_n)$ in $c_0$ instead of $\ell_{\infty}$. Notice that interpolating sequences $(z_n)$ for Bloch type spaces satisfy $|z_n| \rightarrow 1$ since they do not have accumulation points in $\D$. The connection between $c_0-$interpolating sequences and interpolating sequences has been studied in the context of uniform algebras (see \cite{GLM}). In particular, the authors proved that $c_0-$interpolating sequences for $\hinf$ are indeed interpolating for $\hinf$. We will show that this result remains true if we deal with $\Binf$. \medskip


In the proof of the next theorem we will use the following result (see \cite{M}, Theorem 5, p. 82): let $X,Y$ be Banach spaces and $T: X \to Y$ a linear and bounded operator. Then $T$ is bounded below if and only if $T^*$ is surjective. Furthermore, $T$ is surjective if and only if $T^*$ is bounded below.

\begin{teo} \label{teo21}
Let $v$ be a typical weight on $\D$. If $(z_n) \subset \D$ is a sequence of distinct points, then the following statements are equivalent:
\begin{itemize}
\item[(a)] The sequence $(z_n)$ is interpolating for $\Binf$.
\item[(b)] There exists a constant $C > 0$ such that:
$$\|(\xi_n)\|_1 \leq C \left\| \sum_{n=1}^\infty \xi_n v(z_n) \delta'_{z_n}\right\| \ \mbox{ for any } (\xi_n) \in \ell_1.$$
\item[(c)] The sequence $(z_n)$ is $c_0-$interpolating for $\Bo$.
\end{itemize}
\end{teo}

\dem Define $S:\ell_1 \to$ $\Bpre$ given by:
$$S((\xi_n))=\sum_{n=1}^{\infty} \xi_n v(z_n) \delta'_{z_n},$$
\noindent which is clearly a well-defined, linear, continuous map. Condition (b) states that $S$ is bounded below. We have $(\Bpre)^*=\Binf$ by Corollary \ref{rem12} and it is easy that $S^*=T$, where $T$ is the interpolating operator on $\Binf$. \medskip

(a) $\Leftrightarrow$ (b) It is clear since (a) states that the interpolating operator $T$ is surjective and this is equivalent to $S$ being bounded below. \medskip

(b) $\Leftrightarrow$ (c) Notice that $\Bpre=(\Bo)^*$ by Corollary \ref{rem12}. Since $|z_n| \rightarrow 1$, $T$ maps $\Bo$ into $c_0$ and $(T|_{\Bo})^*=S$. Hence (b) and (c) are equivalent by the result above. \qed \bigskip

From Theorem \ref{teo21} and its proof we have the following results:

\begin{cor} \label{coro}
We have:
\begin{itemize}
\item[(a)] A sequence $(z_n) \subset \D$ is interpolating for $\Binf$ if and only if it is $c_0-$interpolating for $\Bo$.
\item[(b)] $(T|_{\Bo})^{**}=T$ and $T$ is $w^*-w^*-$continuous.
\end{itemize}
\end{cor}

The inspiration to the next result comes from Theorem 2.4 in \cite{GLM} and Proposition 7.7 in \cite{BGLM}.

\begin{teo} \label{teo23}
Let $v$ be a typical weight and $(z_n) \subset \D$. Suppose that for any $(a_n) \in c_0$ there exists $f \in \Binf$ such that $v(z_n)f'(z_n)=a_n$ for any $n \in \N$. Then $(z_n)$ is interpolating for $\Binf$.
\end{teo}

\dem Let $a=(a_n) \in \ell_{\infty}$. By Goldstein's theorem, there exists a sequence $\{b^k\} \subset c_0$ such that $b^k \stackrel{w^*}{\rightarrow} a$ when $k \rightarrow \infty$. Consider the interpolating operator $T: \Binf \to \ell_{\infty}$ given by $T(f)=( f'(z_n)v(z_n))$ and take $A=T^{-1}(c_0) \subset \Binf$. Since $A$ is closed in $\Binf$, it follows that $A$ is a Banach space. The linear operator $T|_{A}: A \to c_0$ is surjective by the assumption. The induced operator $\tilde{T}|_A:A/ \ker (T|_A) \to c_0$ is thus bounded, injective and surjective. Therefore by the Open Mapping theorem, there is $M >0$ such that:
$$\inf_{h \in \ker(T|_A)} \|f+h\|_{\Binf} \leq \frac{M}{2} \|b\|_{\infty} \ \mbox{ if } \ T|_A(f)=b \in c_0.$$
\noindent Hence for each $b \in c_0$ there is $f \in A$ such that $T|_A(f)=b$ and $\|f\|_{\Binf}\leq M \|b\|_\infty$. In particular for any $k \in \N$ there exists $g_k \in \Binf$ such that $T(g_k)=b^k$ and $\|g_k\|_{\Binf} \leq M \|b^k\|_{\infty}$. Since $\{b^k\}$ is weak-star convergent, it is bounded in $c_0$, so there is $C >0$ such that $\|g_k\| \leq C$ for all $k \in \N$. Since $\Bpre$ is separable and $(\Bpre)^*=\Binf$, by Alaoglu's theorem there exists a subsequence $(g_{k_m})$ of $(g_k)$ which $w^*-$converges to $g \in \Binf$. By Corollary \ref{coro}, $(T|_{\Bo})^{**}=T$ and $T$ is $w^*-w^*-$continuous, so:
$$a=w^*-\lim_{m \rightarrow \infty} b^{k_m} = w^*-\lim_{m \rightarrow \infty} T(g_{k_m})=T(g).$$ \qed \bigskip

Recall that an idempotent $(a_n)$ of $\ell_\infty$ is a sequence that satisfies $a_n^2=a_n$ for any $n \in \N$, that is, $a_n=0$ or $a_n=1$ for any $n \in \N$. Hayman proved that it is sufficient to interpolate idempotent elements $(a_n) \in \linf$ to assure that the sequence $(z_n) \subset \D$ is interpolating for $\hinf$ (see \cite{Hayman}). We will prove that this result remains true when we deal with interpolating sequences for $\Binf$. \medskip

To prove Theorem \ref{teomain}, we will need the following result due to Beurling (See Theorem 4.3 in \cite{Bade-Curtis}): if $X$ is a Banach space and $L: \ell_1 \to X$ is a linear operator such that $L^*(X)$ is dense in $\linf=(\ell_1)^*$, then $L^*(X^*)=\linf$. \medskip

Now we can state our main result:
\begin{teo} \label{teomain}
Let $v$ be a typical weight and $(z_n)\subset\mathbf{D}$ a sequence of distinct points. Then the following assertions are equivalent:
\begin{itemize}
\item[(a)] $(z_n)$ is interpolating for $\Binf$.
\item[(b)] $(z_n)$ is $c_0-$interpolating for $\Bo$.
\item[(c)] $(z_n)$ is $c_0-$interpolating for $\Binf$.
\item[(d)] $(z_n)$ is interpolating for $\Binf$ when only considering idempotents of $\linf$.
\end{itemize}
\end{teo}

\dem It remains only to prove that $(d) \rightarrow (a)$. Let us consider the interpolating operator  $T:\Binf \to \ell^\infty$. By Theorem \ref{teo21} and Corollary \ref{coro} we have that $(T|_{\Bo})^{**}=T$ and $(T|_{\Bo})^*$ maps $\ell_1$ into $\Bpre$. Therefore we need to prove that $T$ has dense range. Indeed, consider $E\subseteq\mathbb{N}$ and denote by $\chi_E$ the sequence in $\ell^\infty$ given by:
$$\chi_E(n):=\left\{
	\begin{array}{ll}
	1 & \mbox{ if } n\in E\\
	0 & \mbox{ if } n\in \mathbb{N}\setminus E.
	\end{array}\right. $$ 
	
	\noindent Define $S\subset\ell^\infty$ to be the set of functions on $\N$ of the form $\sum_{i=1}^{m}a_i\chi_{A_i}$ such that $m \in \N$, $a_i \in \C$ for any $i=1,2,\ldots,m$ and $(A_i)_{i=1}^m$ are pairwise disjoint sets. The set $S$ is dense in $\ell^\infty$ since $S$ is the set of simple functions in $\ell_\infty=L^\infty(\mathbb{N},c)$ where $c$ is the cardinal measure. Let $x \in S$, that is, $x:=\sum_{i=1}^{m}a_i\chi_{A_i}\in \ell^\infty$. By hypothesis, for any $i=1,2,\ldots,m$ there are functions $f_i\in\Binf$ such that:
$$f_i'(z_n)v(z_n)=\chi_{A_i}(n)\quad \mbox{ for all } n \in \N.$$
\noindent For $f:=\sum_{i=1}^{m}a_i f_i \in\Binf$ we have that:
$$T(f)=(f'(z_n) v(z_n))_{n=1}^\infty=\left( \sum_{i=1}^{m}a_i f_i'(z_n)v(z_n) \right)_{n=1}^\infty=\left( \sum_{i=1}^{m}a_i\chi_{A_i}(n) \right )_{n=1}^\infty=x$$
\noindent and we are done. \qed \bigskip

\subsection*{Examples of interpolating sequences for $\Binf$}

%
Finally we turn to some examples of interpolating sequences for $\Binf$. Recall that a sequence $(z_n) \subset \D$ is said to be a Blaschke sequence if $\sum_{n=1}^{\infty} (1-|z_n|) < \infty$. It is well-known that the sequence of zeros of a non-zero bounded analytic function on $\D$ satisfies the Blaschke condition (see \cite{G}). It is also well-known that interpolating sequences for $\hinf$ satisfy the Blaschke condition. However, there exist interpolating sequences for $\Binf$ for some particular weights which does not satisfy this condition (an easy adaptation of Proposition 6.4 in \cite{BGLM}). 

\begin{prop} \label{add}
If $(z_n) \subset \D$ is interpolating for $\hinf$ and $z_0 \in \D$ is such that $z_0 \neq z_n$ for every $n \in \N$, then the sequence $(w_n)$ defined by $w_1=z_0$ and $w_{n+1}=z_n$ for $n\geq1$  is also interpolating for $\hinf$.
\end{prop}
\dem Let $(z_n) \subset \D$ be an interpolating sequence for $\hinf$. Since it satisfies the Blaschke condition, we can consider its Blaschke product $B: \D \to \C$, which is a bounded analytic function satistfying $B(z)=0$ if and only if $z=z_n$ for some $n \in \N$. Consider $(a_n) \in \linf$. There exists a function $f\in\hinf$ satisfying $f(z_n)=a_{n+1}$ for every $n\in\N$, so if we let $\alpha:=B(z_0)\neq0$ and define $g: \D \to \C$ by:
$$g(z):=\frac{(a_1-f(z_0))}{\alpha} B(z)+f(z)$$
\noindent we have that $g\in\hinf$ and $g(w_n)=a_n$ for all $n\in \N$. \qed \bigskip

Madigan and Matheson proved that if a sequence is sufficiently separated for the pseudohyperbolic distance, then it is interpolating for the classical Bloch space $\B$. We prove that this condition is not necessary since we can find interpolating sequences for $\B_{v_\alpha}^{\infty}$, where $v_{\alpha}(z)=(1-|z|^2)^\alpha$, $\alpha >0$, and in particular for $\B$, as close as we want:
\begin{prop}
Let $\alpha > 0$. For any $\eps >0$ there exist interpolating sequences for $\B_{v_\alpha}^{\infty}$ whose constant of separation is less than $\eps$.
\end{prop}

\dem Let $\eps >0$. Consider an interpolating sequence $(z_n)$ for $\hinf$, for example, $z_n:=1-\frac{1}{2^n}$, and add a point $z_0 \notin (z_n)$ such that $\rho(z_0,z_1)<\eps$. By Proposition \ref{add} the sequence given by $\{z_0\} \cup (z_n)$ is interpolating for $\hinf$ hence interpolating for $\B_{v_\alpha}^{\infty}$ by Theorem 6.3 in \cite{BGLM}. \qed

\end{document}